\documentclass[leqno]{article}
\usepackage{amsfonts}
\usepackage{amsmath,amssymb}
\usepackage{latexsym}
\usepackage{color}

\def \Real{\mbox{\sl I\kern-.166em R}}
\def \Hyp{\mbox{\sl I\kern-.166em H}}
\def \Nat{\mbox{\sl I\kern-.166em N}}
\newtheorem{theorem}{Theorem}[section]
\newtheorem{proposition}[theorem]{Proposition}

\newtheorem{remark}[theorem]{Remark}

\input{tcilatex}

\begin{document}

\title{{\large \textbf{On symmetries of the Heisenberg group}}}
\author{W. Batat and A. Zaeim \thanks{Second author supported by funds of University of Payame-Noor}
}
\date{}
\maketitle

\begin{abstract}
We consider the three-dimensional Heisenberg group, equipped with any left-invariant metric, either Lorentzian or Riemannian. We
completely classify their affine vector fields and investigate their
relationship with Killing vector fields and their casual character. We also classify their Ricci, curvature and matter collineations, proving that there are no nontrivial examples.
\end{abstract}

\bigskip

\textbf{Wafaa Batat}

\textsc{\'{E}}cole Nationale Polytechnique d'Oran

B.P 1523 El M'naouar Oran 31000, Algeria

\textit{E-mail address}: batatwafa@yahoo.fr\medskip


%
%
%
\textbf{Amirhesam Zaeim}

Department of Mathematics, Payame noor University,

P.O. Box 19395-3697, Tehran, Iran

\textit{E-mail address}: zaeim@pnu.ac.ir\medskip



\section{{\protect\large \textbf{Introduction and preliminaries}}}
Study of symmetries, over different geometric spaces is one of interesting topics in geometry and mathematical physics.
As we know, several tensor fields live on a (pseudo-)Riemannian manifold $(M,g)$, each of them codifying some geometric or
physical quantity of the space. A symmetry of a tensor field $T$, is a one-parameter group of diffeomorphisms of $(M,g)$,
which leaves $T$ invariant. By this definition, each symmetry corresponds to a vector field $X$ which satisfies $\mathcal L_XT=0$, where $\mathcal L$
denotes the Lie derivative. Some famous symmetries are: symmetries of the metric tensor $g$ which correspond to the Killing vector fields. Symmetries of the Levi-Civita connection $\nabla$ which correspond to the affine vector fields. Symmetries of the Ricci tensor $\varrho$ which correspond to the Ricci collineations. Symmetries of the curvature tensor $R$ (resp. Weyl conformal tensor $W$) which correspond to the curvature collineations (resp. Weyl collineations) and finally, symmetries of the energy-momentum tensor $T=\varrho-\frac12\tau g$ which correspond to the matter collineations. Matter collineations are more relevant from a physical point of view (see \cite{Cam2,Car}).

Since symmetries are more significant from physical aspects, symmetries were studied on several kinds of space-times (see \cite{C5,Calv,Cam1,C6,C7} and references therein). Symmetry analysis of the Heisenberg groups, which arise in the mathematical description of one-dimensional quantum mechanical systems, are the subject of the present research.

{In this paper, }we consider the left-invariant Riemannian and Lorentzian
metrics admitted by the Heisenberg group $H_{3}$. In Section~2, after
reporting some basic information about Heisenberg group and its
left-invariant metrics in global coordinates, we shall describe the
Levi-Civita connection, the curvature and the Ricci tensor of $H_{3}$. The
study of the affine vector fields of these metrics will be made in the next
Section~3, where we shall also report the classification of the
corresponding Killing vector fields, making some needed corrections. For $
g_{1}$ and $g_{2}$, we shall conclude that affine vector fields are necessarily Killing.
On the other hand, affine vector fields of $(H_{3},g_{3})$ will be shown to
form a twelve-dimensional Lie algebra, including properly the
six-dimensional Lie algebra of its Killing vector fields. {is compatible with the fact that this metric is flat \cite{O'N}.} We shall also
investigate the causal character of Killing and affine vector fields in the
different cases. Then, in Section~4 we shall respectively classify Ricci,
curvature and matter collineations on the three-dimensional Heisenberg group
equipped with Riemannian and Lorentzian left-invariant metrics.

\section{\protect\Large Connection and Curvature of Heisenberg group}

\setcounter{equation}{0}

We recall that the Heisenberg group $H_{3}$ is formed by all real valued
matrices of the form
\begin{equation*}
\left(
\begin{array}{ccc}
1 & x & z \\
0 & 1 & y \\
0 & 0 & 1%
\end{array}%
\right) ,
\end{equation*}%
which is a Lie group diffeomorphic to $%
\mathbb{R}
^{3}$ endowed with the multiplication law%
\begin{equation*}
(x,y,z)(\tilde{x},\tilde{y},\tilde{z})=(x+\tilde{x},y+\tilde{y},z+\tilde{z}-x%
\tilde{y}).
\end{equation*}%
A direct calculation yields that the left-invariant vector fields
corresponding to $(1,0,0),(0,1,0),$ $(0,0,1)$ in the tangent space of the
identity are, respectively $\frac{\partial }{\partial x},\frac{\partial }{%
\partial y}-x\frac{\partial }{\partial z},\frac{\partial }{\partial z}$
where $\left\{ \frac{\partial }{\partial x},\frac{\partial }{\partial y},%
\frac{\partial }{\partial z}\right\} $ are the coordinate vector fields for $%
\mathbb{R}
^{3}$.

In this section we shall write down the Levi-Civita connection, the
curvature and the Ricci tensor of the three-dimensional Heisenberg group
equipped with the left-invariant Riemannian (\cite{Mi}) and Lorentzian
metrics (\cite{R}).

Throughout the paper, we denote by $\nabla $ the Levi-Civita connection of $%
(H_{3},g),$ $R$ its curvature tensor, taken with the sign convention:%
\begin{equation*}
R(X,Y)=\left[ \nabla _{X},\nabla _{Y}\right] -\nabla _{\left[ X,Y\right] },
\end{equation*}%
and by $\varrho $ the Ricci tensor of $(H_{3},g)$ which is defined by
\begin{equation*}
\varrho \left( X,Y\right) =\overset{3}{\underset{k=1}{\sum }}g\left(
e_{k},e_{k}\right) g(R(e_{k},X)Y,e_{k}),
\end{equation*}%
where $\{e_{k}\}_{k=1,2,3}$ is an orthonormal basis.

We also denote by $\tau $ the scalar curvature given by $\tau =\sum
g^{ij}\rho (e_{i},e_{j})$ where $g^{ij}$ are the components of the inverse
matrix of the metric tensor, with respect to the orthonormal basis $%
\{e_{k}\}_{k=1,2,3}.$

\subsection{\protect\large Riemannian setting}

Let $\lambda $ denote a positive real number. Consider the left-invariant
vector fields $e_{1},e_{2},e_{3}$ defined by%
\begin{equation}
e_{1}=\frac{\partial }{\partial y}-x\frac{\partial }{\partial z},\text{ }%
e_{2}=\lambda \frac{\partial }{\partial x},\text{ }e_{3}=\frac{\partial }{%
\partial z},\text{\ }  \label{base0}
\end{equation}%
which is dual to $\left\{ \mathrm{d}y,\frac{1}{\lambda }\mathrm{d}x,x\mathrm{%
d}y+\mathrm{d}z\right\} .$ We have the well known Heisenberg bracket
relations%
\begin{equation*}
\lbrack e_{1},e_{2}]=\lambda e_{3},\ \ [e_{1},e_{3}]=0,\ \ [e_{2},e_{3}]=0.
\end{equation*}%
We now endow $H_{3}$ with the Riemannian metric%
\begin{equation}
g_{0}=\frac{1}{\lambda ^{2}}\mathrm{d}x^{2}+\mathrm{d}y^{2}+(x\mathrm{d}y+%
\mathrm{d}z)^{2}.  \label{g0}
\end{equation}%
The metric $g_{0}$ is invariant with respect to the left-translations
corresponding to that multiplication law. Any left-invariant Riemannian metric on the Heisenberg group is isometric to the above left-invariant metric $g_0$. To note that, the left-invariant vector fields given in \eqref{base0} are
orthonormal with respect to $g_{0}.$

The only non-vanishing components of the Levi-Civita connection of $%
(H_{3},g_{0})$ are the following ones:%
\begin{equation}
\nabla _{e_{1}}e_{2} =-\nabla _{e_{2}}e_{1}=\frac{\lambda }{2}e_{3},\quad
\nabla _{e_{1}}e_{3} =\nabla _{e_{3}}e_{1}=-\frac{\lambda }{2}e_{2},\quad
\nabla _{e_{2}}e_{3} =\nabla _{e_{3}}e_{2}=\frac{\lambda }{2}e_{1}.
\label{levi0}
\end{equation}%
The non-vanishing curvature components are:

\begin{equation}
\begin{array}{l}
R(e_{1},e_{2})e_{1} =\frac{3\lambda ^{2}}{4}e_{2},\quad
R(e_{1},e_{2})e_{2} =-\frac{3\lambda ^{2}}{4}e_{1}, \quad
R(e_{1},e_{3})e_{1} =R(e_{2},e_{3})e_{2}=-\frac{\lambda ^{2}}{4}e_{3}, \\
R(e_{1},e_{3})e_{3} =\frac{\lambda ^{2}}{4}e_{1}, \quad
R(e_{2},e_{3})e_{3} =\frac{\lambda ^{2}}{4}e_{2},
\end{array}
\label{Rieman0}
\end{equation}%
and the ones obtained from them using the symmetries of the curvature
tensor. By contraction, the non-zero components $\varrho_{ij}=\varrho (e_{i},e_{j})$ of the
Ricci tensor are:%
\begin{equation}
\varrho _{11}=\varrho _{22}=-\frac{\lambda ^{2}}{2},\quad \varrho _{33}=%
\frac{\lambda ^{2}}{2}.\label{ro0}
\end{equation}%
Consequently, the scalar curvature of $g_{0}$ is given by \
\begin{equation}
\tau =-\frac{\lambda ^{2}}{2}.  \label{to0}
\end{equation}

\subsection{\protect\large Lorentzian setting}

As Rahmani proved in \cite{R}, up to isomorphism on its Lie algebra, $H_{3}$
admits three different left-invariant Lorentzian metrics, to which we shall
refer as $g_{1}$, $g_{2}$ and $g_{3}$ (the flat one). Explicitly, these
metrics are given by:
\begin{eqnarray}
&&g_{1}=-\frac{1}{\lambda ^{2}}\mathrm{d}x^{2}+\mathrm{d}y^{2}+(x\mathrm{d}y+%
\mathrm{d}z)^{2},\quad \lambda \neq 0,  \label{g1} \\
&&g_{2}=\frac{1}{\lambda ^{2}}\mathrm{d}x^{2}+\mathrm{d}y^{2}-(x\mathrm{d}y+%
\mathrm{d}z)^{2},\quad \ \ \lambda \neq 0,  \label{g2} \\
&&g_{3}=\mathrm{d}x^{2}+(x\mathrm{d}y+\mathrm{d}z)^{2}-((1-x)\mathrm{d}y-%
\mathrm{d}z)^{2}.  \label{g3}
\end{eqnarray}%
This is a strong contrast to the Riemannian case in which there is only one
(up to positive homothety) and it is not flat.

\subsubsection{{\protect\large The metric }$g_{1}$}

First, we endow $H_{3}$ with the left-invariant Lorentzian metric $g_{1}$.
The Lie algebra of $H_{3}$ has an orthonormal basis
\begin{equation}
e_{1}=\frac{\partial }{\partial z},\quad e_{2}=\frac{\partial }{\partial y}-x%
\frac{\partial }{\partial z},\quad e_{3}=\lambda \frac{\partial }{\partial x}%
,  \label{base1}
\end{equation}%
with%
\begin{equation*}
g_{1}(e_{1},e_{1})=g_{1}(e_{2},e_{2})=1,\text{ \ \ }g_{1}(e_{3},e_{3})=-1.
\end{equation*}%
for which, we have the Lie products:
\begin{equation}
\lbrack e_{2},e_{3}]=\lambda e_{1},\ \ [e_{1},e_{2}]=0,\ \ [e_{1},e_{3}]=0.
\label{crochet1}
\end{equation}%
The only non-vanishing components of the Levi-Civita connection of $%
(H_{3},g_{1})$ are as follows:%
\begin{equation}
\nabla _{e_{1}}e_{2} =\nabla _{e_{2}}e_{1}=\frac{\lambda }{2}e_{3}, \quad
\nabla _{e_{1}}e_{3} =\nabla _{e_{3}}e_{1}=\frac{\lambda }{2}e_{2}, \quad
\nabla _{e_{2}}e_{3} =-\nabla _{e_{3}}e_{2}=\frac{\lambda }{2}e_{1}.
\label{levi1}
\end{equation}%
The only non-vanishing curvature components are
\begin{equation}
\begin{array}{l}
R(e_{1},e_{2})e_{1} =\frac{\lambda ^{2}}{4}e_{2}, \quad
R(e_{1},e_{2})e_{2} =-R(e_{1},e_{3})e_{3}=-\frac{\lambda ^{2}}{4}e_{1}, \quad
R(e_{1},e_{3})e_{1} =\frac{\lambda ^{2}}{4}e_{3}, \\
R(e_{2},e_{3})e_{2} =-\frac{3\lambda ^{2}}{4}e_{3}, \quad
R(e_{2},e_{3})e_{3} =-\frac{3\lambda ^{2}}{4}e_{2},
\end{array}
\label{Rieman1}
\end{equation}%
and the ones obtained from them using the symmetries of the curvature
tensor. Consequently, the components $\varrho _{ij}=\varrho (e_{i},e_{j})$
of the Ricci tensor are:%
\begin{equation}
\varrho _{11}=\varrho _{33}=-\frac{\lambda ^{2}}{2},\quad \varrho _{22}=%
\frac{\lambda ^{2}}{2},\quad \varrho _{ij}=0\text{ otherwise.}  \label{ro1}
\end{equation}%
The scalar curvature of $g_{1}$ is given by \
\begin{equation}
\tau =\frac{\lambda ^{2}}{2}.  \label{to1}
\end{equation}

\subsubsection{{\protect\large The metric }$g_{2}$}

Now we endow $H_{3}$ with the left-invariant Lorentzian metric $g_{2}$.
Consider the left-invariant orthonormal basis%
\begin{equation}
e_{1}=\frac{\partial }{\partial y}-x\frac{\partial }{\partial z},\quad
e_{2}=\lambda \frac{\partial }{\partial x},\quad e_{3}=\frac{\partial }{%
\partial z},  \label{base2}
\end{equation}%
for which
\begin{equation}
\lbrack e_{1},e_{2}]=\lambda e_{3},\quad \lbrack e_{1},e_{3}]=0,\quad
\lbrack e_{2},e_{3}]=0,  \label{crochet2}
\end{equation}%
and%
\begin{equation*}
g_{2}(e_{1},e_{1})=g_{2}(e_{2},e_{2})=1,\quad g_{2}(e_{3},e_{3})=-1.
\end{equation*}%
The only non-vanishing components of the Levi-Civita connection $\nabla $,
the Riemann curvature tensor $R$ and the Ricci tensor $\varrho $ of $%
(H_{3},g_{2}),$ are respectively given by%
\begin{equation}
\nabla _{e_{1}}e_{2} =-\nabla _{e_{2}}e_{1}=\frac{\lambda }{2}e_{3}, \quad
\nabla _{e_{1}}e_{3} =\nabla _{e_{3}}e_{1}=\frac{\lambda }{2}e_{2}, \quad
\nabla _{e_{2}}e_{3} =\nabla _{e_{3}}e_{2}=-\frac{\lambda }{2}e_{1},
\label{levi2}
\end{equation}

\begin{equation}
\begin{array}{l}
R(e_{1},e_{2})e_{1} =-\frac{3\lambda ^{2}}{4}e_{2}, \quad
R(e_{1},e_{2})e_{2} =\frac{3\lambda ^{2}}{4}e_{1}, \quad
R(e_{1},e_{3})e_{1} =R(e_{2},e_{3})e_{2}=\frac{\lambda ^{2}}{4}e_{3}, \\
R(e_{1},e_{3})e_{3} =\frac{\lambda ^{2}}{4}e_{1}, \quad
R(e_{2},e_{3})e_{3} =\frac{\lambda ^{2}}{4}e_{2},
\end{array}
\label{Rieman2}
\end{equation}%
and%
\begin{equation}
\varrho _{11}=\varrho _{22}=\varrho _{33}=\frac{\lambda ^{2}}{2},\quad
\varrho _{ij}=0\text{ otherwise.}  \label{ro2}
\end{equation}%
The scalar curvature of $g_{2}$ is determined by \
\begin{equation}
\tau =\frac{\lambda ^{2}}{2}.  \label{to2}
\end{equation}

\subsubsection{{\protect\large The metric }$g_{3}$}

Now we consider $H_{3}$ equipped with the left-invariant Lorentzian metric $%
g_{3}$. The Lie algebra of $H_{3}$ has an orthonormal basis%
\begin{equation}
e_{1}=\frac{\partial }{\partial x},\ \ \ \ e_{2}=\frac{\partial }{\partial y}%
+(1-x)\frac{\partial }{\partial z},\ \ \ e_{3}=\frac{\partial }{\partial y}-x%
\frac{\partial }{\partial z},  \label{base3}
\end{equation}%
with%
\begin{equation*}
g_{3}(e_{1},e_{1})=g_{3}(e_{2},e_{2})=1,\quad g_{3}(e_{3},e_{3})=-1,
\end{equation*}%
for which
\begin{equation}
\lbrack e_{2},e_{3}]=0,\ \ \ [e_{2},e_{1}]=[e_{3},e_{1}]=e_{2}-e_{3}.
\label{crochet3}
\end{equation}%
The only non-vanishing components of the Levi-Civita connection $\nabla $ are
given by
\begin{equation}
\nabla _{e_{2}}e_{1} =\nabla _{e_{3}}e_{1}=e_{2}-e_{3},\quad
\nabla _{e_{2}}e_{2} =\nabla _{e_{2}}e_{3}=\nabla _{e_{3}}e_{3}=\nabla
_{e_{3}}e_{2}=-e_{1}.  \notag
\label{levi3}
\end{equation}
The Riemann curvature tensor $R$ (and so, the Ricci tensor $\varrho $) of $%
(H_{3},g_{3})$ vanishes identically.

\section{{\protect\large \textbf{Affine vector fields}}}

\setcounter{equation}{0}

Given a pseudo-Riemannian manifold $(M,g)$, a vector field $X$ tangent to $M$
is Killing if it satisfies one of the following equivalent conditions:

\begin{itemize}
\item[(i)] its local fluxes are given by isometries,

\item[(ii)] $\mathcal{L }_X g=0$, where $\mathcal{L}$ denotes the Lie
derivative,

\item[(iii)] for all vector fields $Y,Z$ tangent to $M$:
\begin{equation}  \label{Kill}
g(\nabla _Y X,Z)=-g(\nabla _Z X,Y).
\end{equation}
\end{itemize}

Killing vector fields have been extensively studied, both in Riemannian and
in Lorentzian geometries, under several different points of view. Some of
the many relevant topics where Killing vector fields played an important
role, are given by functional analysis, Lie groups actions and Bochner's
technique. Summarizing results concerning Killing vector fields would be too
big a task for this paper. We may refer to \cite{S} for a good survey
in the Lorentzian case. It is worthwhile to emphasize that for a Lorentzian
manifold, the causal character of a Killing vector field aquires a special
relevance, and the most important causality condition is given by a timelike
Killing vector field.

Next, a vector field $X$ tangent to $(M,g)$ is said to be \emph{affine} if
it satisfies one of the following equivalent conditions:

\begin{itemize}
\item[(i)] its local fluxes are given by affine maps,

\item[(ii)] $\mathcal{L }_X \nabla=0$, where $\nabla$ is the Levi Civita
connection of $(M,g)$,

\item[(iii)] for all vector fields $Y,Z$ tangent to $M$:
\begin{equation}  \label{aff}
[X,\nabla _Y Z] =\nabla _{[X,Y]} Z + \nabla _Y [X,Z].
\end{equation}
\end{itemize}

Clearly, a Killing vector field is also affine and it is an interesting
problem whether the converse holds for a given manifold $(M,g)$. In
particular, if $(M,g)$ is a simply connected spacetime, the existence of a
\emph{proper} (that is, non Killing) affine vector field implies the
existence of a second-order covariantly constant symmetric tensor, nowhere
vanishing, not proportional to $g$. As a consequence, the holonomy group of
the manifold is reducible (see for example \cite{Sh}).

In this Section, we completely classify affine vector fields for the
left-invariant Riemannian and Lorentzian metrics admitted by the Heisenberg
group $H_{3}$.

\subsection{\protect\large Riemannian setting}

The Lie algebra of Killing vector field of $(H_{3},g_{0})$ is generated by
the following vectors
\begin{equation}
X_{1}=\lambda ^{2}y\partial _{x}-x\partial _{y}-\frac{1}{2}(\lambda
^{2}y^{2}-x^{2})\partial _{z},\quad X_{2}=\partial _{x}-y\partial _{z},\quad
X_{3}=\partial _{y},\quad X_{4}=\partial _{z}.  \label{g0kill}
\end{equation}

Next, {let }$X=f_{1}e_{1}+f_{2}e_{2}+f_{3}e_{3}$ be a tangent vector field
on $(H_{3},g_{0}),$ where $g_{0}$ is described by \eqref{g0} and $%
f_{i}=f_{i}(x,y,z),i=1,2,3$ are arbitrary smooth functions. We apply %
\eqref{aff} when $Y,Z$ are chosen in the {orthonormal basis \eqref{base0}}.
Via standard calculations, we find that $X$ is an affine vector field of $%
(H_{3},g_{1})$ if and only if $f_{1},f_{2},f_{3}$ satisfy the following
system of partial differential equations:%
\begin{equation}
\begin{array}{l}
\frac{\partial ^{2}f_{1}}{\partial y^{2}}-2x\frac{\partial ^{2}f_{1}}{%
\partial y\partial z}+x^{2}\frac{\partial ^{2}f_{1}}{\partial z^{2}}=0, \\
\frac{\partial ^{2}f_{2}}{\partial y^{2}}-2x\frac{\partial ^{2}f_{2}}{%
\partial y\partial z}+x^{2}\frac{\partial ^{2}f_{2}}{\partial z^{2}}-\lambda
\left( \lambda f_{2}+\frac{\partial f_{3}}{\partial y}-x\frac{\partial f_{3}%
}{\partial z}\right) =0, \\
\frac{\partial ^{2}f_{3}}{\partial y^{2}}-2x\frac{\partial ^{2}f_{3}}{%
\partial y\partial z}+\lambda \left( \frac{\partial f_{2}}{\partial y}-x%
\frac{\partial f_{2}}{\partial z}\right) =0, \\
2\left( \frac{\partial ^{2}f_{1}}{\partial x\partial y}-x\frac{\partial
^{2}f_{1}}{\partial x\partial z}\right) +\lambda f_{2}+\frac{\partial f_{3}}{%
\partial y}-x\frac{\partial f_{3}}{\partial z}-\frac{\partial f_{1}}{%
\partial z}=0, \\
2\left( \frac{\partial ^{2}f_{2}}{\partial x\partial y}-x\frac{\partial
^{2}f_{2}}{\partial x\partial z}\right) +\lambda f_{1}-\lambda \frac{%
\partial f_{3}}{\partial x}-\frac{\partial f_{2}}{\partial z}=0, \\
2\left( \frac{\partial ^{2}f_{3}}{\partial x\partial y}-x\frac{\partial
^{2}f_{3}}{\partial x\partial z}\right) +x\frac{\partial f_{1}}{\partial z}-%
\frac{\partial f_{1}}{\partial y}+\lambda \frac{\partial f_{2}}{\partial x}-%
\frac{\partial f_{3}}{\partial z}=0, \\
2\left( \frac{\partial ^{2}f_{1}}{\partial y\partial z}-x\frac{\partial
^{2}f_{1}}{\partial z^{2}}\right) +\lambda \left( \frac{\partial f_{2}}{%
\partial y}-x\frac{\partial f_{2}}{\partial z}+\lambda \frac{\partial f_{1}}{%
\partial x}\right) =0, \\
2\left( \frac{\partial ^{2}f_{2}}{\partial y\partial z}-x\frac{\partial
^{2}f_{2}}{\partial z^{2}}\right) +\lambda \left( -\frac{\partial f_{3}}{%
\partial z}+x\frac{\partial f_{1}}{\partial z}-\frac{\partial f_{1}}{%
\partial y}+\lambda \frac{\partial f_{2}}{\partial x}\right) =0, \\
2\frac{\partial ^{2}f_{3}}{\partial y\partial z}+\lambda \left( \frac{%
\partial f_{2}}{\partial z}-\lambda f_{1}+\lambda \frac{\partial f_{3}}{%
\partial x}\right) =0, \\
\frac{\partial f_{3}}{\partial x}+\frac{\partial ^{2}f_{1}}{\partial x^{2}}%
-f_{1}=0, \\
\frac{\partial ^{2}f_{2}}{\partial x^{2}}=0, \\
\frac{\partial ^{2}f_{3}}{\partial x^{2}}-\frac{\partial f_{1}}{\partial x}%
=0, \\
\frac{\partial f_{3}}{\partial z}+\lambda \frac{\partial f_{2}}{\partial x}+x%
\frac{\partial f_{1}}{\partial z}-\frac{\partial f_{1}}{\partial y}+2\frac{%
\partial f_{1}}{\partial x\partial z}=0, \\
\lambda \frac{\partial f_{1}}{\partial x}-x\frac{\partial f_{2}}{\partial z}+%
\frac{\partial f_{2}}{\partial y}-2\frac{\partial ^{2}f_{2}}{\partial
x\partial z}=0, \\
\lambda f_{2}+\frac{\partial f_{3}}{\partial y}-x\frac{\partial f_{3}}{%
\partial z}+\frac{\partial f_{1}}{\partial z}-2\frac{\partial ^{2}f_{3}}{%
\partial x\partial z}=0, \\
\lambda \frac{\partial f_{2}}{\partial z}+\frac{\partial ^{2}f_{1}}{\partial
z^{2}}=0, \\
\lambda \frac{\partial f_{1}}{\partial z}-\frac{\partial ^{2}f_{2}}{\partial
z^{2}}=0, \\
\frac{\partial ^{2}f_{3}}{\partial z^{2}}=0.%
\end{array}
\label{affg0}
\end{equation}%
{
By straightforward calculations, using the above system of PDEs, we get
$$
\begin{array}{l}
\frac{\partial ^{2}f_{1}}{\partial x^2}=\frac{\partial f_{1}}{\partial y}=\frac{\partial f_{1}}{\partial z}=0,
\end{array}
$$
so by integrating we have}
\begin{equation*}
f_{1}=c_{1}x+c_{2},\text{ }c_{1},c_{2}\in
\mathbb{R}
.
\end{equation*}%
%
{
Now, by direct computations, we obtain  $\frac{%
\partial f_{1}}{\partial z}=\frac{\partial f_{2}}{\partial z}=0,$ $\frac{%
\partial f_{2}}{\partial y}=-\lambda \frac{\partial f_{1}}{\partial x},$ which gives (after replacing the expression of $f_{1}$)
}
\begin{equation*}
f_{2}=-\lambda c_{1}y+c_{3},\text{ }c_{3}\in
\mathbb{R}
.
\end{equation*}%
By replacing $f_{1}$ and $f_{2}$ into the fourth equation of \eqref{affg0}
we prove \ that
\begin{equation*}
f_{3}=\frac{\lambda ^{2}c_{1}}{2}y^{2}-\lambda c_{3}y+A(x),
\end{equation*}%
for a smooth function $A=A(x)$. Hence, { since have $f_1=\frac{\partial f_3}{\partial x}$, }
%
it implies that $%
A^{\prime }(x)=c_{1}x+c_{2},$ and so,
\begin{equation*}
f_{3}=\frac{\lambda ^{2}c_{1}}{2}y^{2}+\frac{c_{1}}{2}x^{2}+c_{2}x-\lambda
c_{3}y+c_{4},\text{ }c_{4}\in
\mathbb{R}
.
\end{equation*}

So, we proved the following

\begin{theorem}
$X=f_{1}e_{1}+f_{2}e_{2}+f_{3}e_{3}$ is an affine vector field on $%
(H_{3},g_{0})$ if and only if
\begin{equation}
\left\{
\begin{array}{l}
f_{1}(x,y,z)=c_{1}x+c_{2}, \\
f_{2}(x,y,z)=-\lambda c_{1}y+c_{3}, \\
f_{3}(x,y,z)=\frac{\lambda ^{2}c_{1}}{2}y^{2}+\frac{c_{1}}{2}%
x^{2}+c_{2}x-\lambda c_{3}y+c_{4},%
\end{array}%
\right.   \label{affg0com}
\end{equation}%
for any real constants $c_{1},..,c_{4}.$
\end{theorem}

\noindent {Next, one can check easily that the Lie algebra of affine vector
fields of }$(H_{3},g_{0}),$ admits as basis the vector fields given in %
\eqref{g0kill}. Therefore, the following result holds:

\begin{theorem}
Killing vector fields are all and the ones affine vector fields on $%
(H_{3},g_{0})$.
\end{theorem}

\subsection{\protect\Large Lorentzian setting}

\subsubsection{{\protect\large \textbf{The first metric}}}

{Killing vector fields of }$(H_{3},g_{1})$ were already classified in \cite%
{R}. More precisely, it was shown that the Lie algebra of Killing vector
fields of $(H_{3},g_{1})$ is four-dimensional and admits the following
basis:
\begin{equation}
X_{1}=\lambda ^{2}y\partial _{x}+x\partial _{y}-\frac{1}{2}(x^{2}+\lambda
^{2}y^{2})\partial _{z},\quad X_{2}=\partial _{x}-y\partial _{z},\quad
X_{3}=\partial _{y},\quad X_{4}=\partial _{z}.  \label{g1kill}
\end{equation}%
{Notice that}, $X_{3},X_{4}$ are spacelike vector fields, while $%
||X_{2}||^{2}=y^{2}-\frac{1}{\lambda ^{2}}$ and $||X_{1}||^{2}=(x^{2}-%
\lambda ^{2}y^{2})(1+\frac{x^{2}-\lambda ^{2}y^{2}}{4})$ and so, the causal
character of $X_{1},X_{2}$ depends on the point. It is easy to check that
neither timelike nor null vector fields can be obtained as linear
combination of $X_{1},X_{2},X_{3},X_{4}$. Hence, we have the following

\begin{proposition}
The Lie algebra of Killing vector fields of $(H_3,g_1)$ contains spacelike
vector fields and vector fields with changing causal character, but neither
timelike nor null vector fields.
\end{proposition}

Next, {let }$X=f_{1}e_{1}+f_{2}e_{2}+f_{3}e_{3}$ be a tangent vector field
on $(H_{3},g_{1}),$ where $f_{i}=f_{i}(x,y,z),i=1,2,3$ are arbitrary smooth
functions. We apply \eqref{aff} when $Y,Z$ are chosen in the {orthonormal
basis \eqref{base1}}. Via standard calculations we find that $X$ is an
affine vector field of $(H_{3},g_{1})$ if and only if $f_{1},f_{2},f_{3}$
satisfy the following system of eighteen partial differential equations:

\begin{equation}
\begin{array}{l}
\frac{\lambda ^{2}}{2}(-\frac{\partial f_{1}}{\partial x}+f_{2})+\frac{%
\lambda }{2}(\frac{\partial f_{3}}{\partial z})+\frac{\partial }{\partial z}(%
\frac{\partial f_{1}}{\partial y}-x\frac{\partial f_{1}}{\partial z})=0, \\
\frac{\partial }{\partial z}(\frac{\partial f_{2}}{\partial y}-x\frac{%
\partial f_{2}}{\partial z})+\frac{\lambda }{2}(\frac{\partial f_{3}}{%
\partial y}-x\frac{\partial f_{3}}{\partial z})-\frac{\lambda ^{2}}{2}(\frac{%
\partial f_{2}}{\partial x})=0 \\
\frac{\lambda }{2}(\frac{\partial f_{1}}{\partial z})-\frac{\lambda ^{2}}{2}(%
\frac{\partial f_{3}}{\partial x})+\frac{\partial }{\partial z}(\frac{%
\partial f_{3}}{\partial y}-x\frac{\partial f_{3}}{\partial z})+\frac{%
\lambda }{2}(\frac{\partial f_{2}}{\partial y}-x\frac{\partial f_{2}}{%
\partial z})=0 \\
x\frac{\partial f_{1}}{\partial z}-\frac{\partial f_{1}}{\partial y}-\lambda
f_{3}+2\frac{\partial ^{2}f_{1}}{\partial z\partial x}-\frac{\partial f_{2}}{%
\partial z}=0 \\
x\frac{\partial f_{2}}{\partial z}-\frac{\partial f_{2}}{\partial y}+\frac{%
\partial f_{1}}{\partial z}+2\frac{\partial ^{2}f_{2}}{\partial{z}\partial
x}+\lambda \frac{\partial f_{3}}{\partial x}=0 \\
x\frac{\partial f_{3}}{\partial z}-\frac{\partial f_{3}}{\partial y}+2\frac{%
\partial ^{2}f_{3}}{\partial z\partial x}+\lambda \frac{\partial f_{2}}{%
\partial x}=0 \\
x\frac{\partial f_{2}}{\partial z}-\frac{\partial f_{1}}{\partial z}-\frac{%
\partial f_{2}}{\partial y}+2(\frac{\partial ^{2}f_{1}}{\partial y\partial x}%
-x\frac{\partial ^{2}f_{1}}{\partial z\partial x})+\lambda \frac{\partial
f_{3}}{\partial x}=0 \\
\frac{\partial f_{1}}{\partial y}-\frac{\partial f_{2}}{\partial z}-x\frac{%
\partial f_{1}}{\partial z}+\lambda f_{3}+2(\frac{\partial ^{2}f_{2}}{%
\partial y\partial x}-x\frac{\partial ^{2}f_{2}}{\partial z\partial x})=0 \\
\lambda (\frac{\partial f_{1}}{\partial x}-f_{2})+2(\frac{\partial ^{2}f_{3}%
}{\partial y\partial x}-x\frac{\partial ^{2}f_{3}}{\partial z\partial x})-%
\frac{\partial f_{3}}{\partial z}=0 \\
\frac{\partial ^{2}f_{1}}{\partial z^{2}}=0 \\
\frac{\partial ^{2}f_{2}}{\partial y^{2}}-2x\frac{\partial ^{2}f_{2}}{%
\partial y\partial z}+x^{2}\frac{\partial ^{2}f_{2}}{\partial z^{2}}=0 \\
\frac{\partial ^{2}f_{3}}{\partial x^{2}}=0 \\
\frac{\partial ^{2}f_{2}}{\partial z^{2}}+\lambda \frac{\partial f_{3}}{%
\partial z}=0 \\
\frac{\partial ^{2}f_{3}}{\partial z^{2}}+\lambda \frac{\partial f_{2}}{%
\partial z}=0 \\
\frac{\partial ^{2}f_{1}}{\partial y^{2}}-2x\frac{\partial ^{2}f_{1}}{%
\partial y\partial z}+x^{2}\frac{\partial ^{2}f_{1}}{\partial z^{2}}+\lambda
(\frac{\partial f_{3}}{\partial y}-x\frac{\partial f_{3}}{\partial z})=0 \\
\lambda (\frac{\partial f_{1}}{\partial y}-x\frac{\partial f_{1}}{\partial z}%
+\lambda f_{3})+\frac{\partial ^{2}f_{3}}{\partial y^{2}}-2x\frac{\partial
^{2}f_{3}}{\partial y\partial z}+x^{2}\frac{\partial ^{2}f_{3}}{\partial
z^{2}}=0 \\
\frac{\partial ^{2}f_{1}}{\partial x^{2}}-\frac{\partial f_{2}}{\partial x}=0
\\
\frac{\partial f_{1}}{\partial x}+\frac{\partial ^{2}f_{2}}{\partial x^{2}}%
-f_{2}=0%
\end{array}
\label{affg1}
\end{equation}%
{We will completely solve the system of PDE }\eqref{affg1}, determining the
affine vector fields of the three-dimensional Heisenberg group $%
(H_{3},g_{1}).$

{By direct computations, using the above system of PDEs, we may find
$\frac{\partial^3 f_1}{\partial y^3}=\frac{\partial^2 f_1}{\partial x\partial y}=\frac{\partial f_1}{\partial z}=0$,
 so}
we have $\frac{\partial f_{1}}{\partial y}=c_{1}y+c_{2}$, $c_{1},c_{2}\in\mathbb{R}.$ Thus,%
\begin{equation*}
f_{1}=\frac{c_{1}}{2}y^{2}+c_{2}y+A(x),
\end{equation*}%
where $A$ is a smooth function of $x.$

{ Taking the derivative of 
$f_2=\frac{\partial f_1}{\partial x}$} with respect to $x$ and using {$f_3=-\frac{1}{\lambda}\frac{\partial f_1}{\partial y}$ and $\frac{\partial f_3}{\partial y}=\lambda\frac{\partial f_2}{\partial x}$}
, we obtain $\frac{\partial ^{2}f_{1}}{\partial x^2}=-\frac{1}{\lambda ^{2}}\frac{\partial ^{2}f_{1}}{\partial y^2}.$ Using the
expression of $f_{1}$ we prove that $A(x)=-\frac{c_{1}}{2\lambda ^{2}}%
x^{2}+c_{3}x+c_{4},$ $c_{3},c_{4}\in
\mathbb{R}
.$
Replacing into {$f_2=\frac{\partial f_1}{\partial x}$, $f_3=-\frac{1}{\lambda}\frac{\partial f_1}{\partial y}$} 
, we then conclude%
\begin{equation*}
f_{2}=-\frac{c_{1}}{\lambda ^{2}}x+c_{3}\text{ \ and \ }f_{3}=-\frac{1}{%
\lambda }\left( c_{1}y+c_{2}\right) .
\end{equation*}

In this way, we proved the following

\begin{theorem}
\label{thg1} $X=f_{1}e_{1}+f_{2}e_{2}+f_{3}e_{3}$ is an affine vector field
on $(H_{3},g_{1})$ if and only if
\begin{equation}
\left\{
\begin{array}{l}
f_{1}(x,y,z)=\frac{c_{1}}{2}y^{2}-\frac{c_{1}}{2\lambda ^{2}}%
x^{2}+c_{3}x+c_{2}y+c_{4}, \\
f_{2}(x,y,z)=-\frac{c_{1}}{\lambda ^{2}}x+c_{3}, \\
f_{3}(x,y,z)=-\frac{1}{\lambda }\left( c_{1}y+c_{2}\right) ,%
\end{array}%
\right.  \label{affg1comp}
\end{equation}%
for any real constants $c_{1},..,c_{4}$.
\end{theorem}

\noindent {A straightforward computations leads to conclude that the Lie algebra
of affine vector fields of }$(H_{3},g_{1}),$ admits a basis of affine vector
fields which coincide with vectors given in \eqref{g1kill}. Therefore, the
following result holds:

\begin{theorem}
\label{thg1a1} Killing vector fields are all and the ones affine vector
fields on $(H_{3},g_{1})$.
\end{theorem}

\subsubsection{{\protect\large \textbf{The second metric}}}

\setcounter{equation}{0}

The four-dimensional Lie algebra of Killing vector fields on $(H_{3},g_{2})$
was determined in \cite{R}. With respect to coordinate vector fields $%
\partial _{x},\partial _{y},\partial _{z}$, one basis of this Lie algebra is
described as follows:
\begin{equation}
X_{1}=-\lambda ^{2}y\partial _{x}+x\partial _{y}+\frac{1}{2}(\lambda
^{2}y^{2}-x^{2})\partial _{z},\quad X_{2}=\partial _{x}-y\partial _{z},\quad
X_{3}=\partial _{y},\quad X_{4}=\partial _{z}.  \label{g2kill}
\end{equation}%
Notwithstanding the formal similarity and many common properties between
Lorentzian metrics $g_{1}$ and $g_{2}$ on $H_{3}$, the causal character of
their Killing vector fields show a relevant difference between them. In
fact, now $X_{4}$ is a timelike vector field, while $||X_{1}||^{2}=(x^{2}+%
\lambda ^{2}y^{2})(1-\frac{x^{2}+\lambda ^{2}y^{2}}{4})$, $||X_{2}||^{2}=%
\frac{1}{\lambda ^{2}}-y^{2}$ and $||X_{3}||^{2}=1-x^{2}$. Thus, the causal
character of $X_{1},X_{2},X_{3}$ depends on the point. One can check that
neither spacelike nor null vector fields can be obtained as linear
combination of $X_{1},X_{2},X_{3},X_{4}$. Therefore, we have the following

\begin{proposition}
The Lie algebra of Killing vector fields of $(H_3,g_2)$ contains timelike
vector fields and vector fields with changing causal character, but neither
spacelike nor null vector fields.
\end{proposition}

We now determine affine vector fields of $(H_{3},g_{2})$. Apply \eqref{aff}
when $Y,Z$ are chosen in the orthonormal basis $\{e_{1},e_{2},e_{3}\}$ given
by \eqref{base2}, we find that $X=f_{1}e_{1}+f_{2}e_{2}+f_{3}e_{3},$ where $%
f_{i}=f_{i}(x,y,z),i=1,2,3$ are arbitrary smooth functions, is an affine
vector field of $(H_{3},g_{2})$ if and only if $f_{1},f_{2},f_{3}$ satisfy
the following system of eighteen partial differential equations:
\begin{equation}
\begin{array}{l}
\frac{\partial ^{2}f_{1}}{\partial y^2}-2x\frac{\partial ^{2}f_{1}}{%
\partial y\partial z}+x^{2}\frac{\partial ^{2}f_{1}}{\partial z^2}=0, \\
\frac{\partial ^{2}f_{2}}{\partial y^2}-2x\frac{\partial ^{2}f_{2}}{%
\partial y\partial z}+x^{2}\frac{\partial ^{2}f_{2}}{\partial z^2}+\lambda
(\lambda f_{2}+\frac{\partial f_{3}}{\partial y}-x\frac{\partial f_{3}}{%
\partial z})=0, \\
\frac{\partial ^{2}f_{3}}{\partial y^2}-2x\frac{\partial ^{2}f_{3}}{%
\partial y\partial z}+x^{2}\frac{\partial ^{2}f_{3}}{\partial z^2}+\lambda
(\frac{\partial f_{2}}{\partial y}-x\frac{\partial f_{2}}{\partial z})=0, \\
\frac{\partial ^{2}f_{2}}{\partial x^2}=0, \\
\frac{\partial ^{2}f_{1}}{\partial x^2}-\frac{\partial f_{3}}{\partial x}%
+f_{1}=0, \\
\frac{\partial ^{2}f_{3}}{\partial x^2}-\frac{\partial f_{1}}{\partial x}%
=0, \\
\frac{\partial ^{2}f_{3}}{\partial z^2}=0, \\
\frac{\partial ^{2}f_{1}}{\partial z^2}-\lambda \frac{\partial f_{2}}{%
\partial z}=0, \\
\frac{\partial ^{2}f_{2}}{\partial z^2}+\lambda \frac{\partial f_{1}}{%
\partial z}=0, \\
\frac{\partial f_{1}}{\partial z}+\frac{\partial f_{3}}{\partial y}-x\frac{%
\partial f_{3}}{\partial z}+\lambda f_{2}-2(\frac{\partial ^{2}f_{1}}{%
\partial y\partial x}-x\frac{\partial ^{2}f_{1}}{\partial x\partial z})=0,
\\
\frac{\partial f_{2}}{\partial z}+\lambda (f_{1}-\frac{\partial f_{3}}{%
\partial x})-2(\frac{\partial ^{2}f_{2}}{\partial x\partial y}-x\frac{%
\partial ^{2}f_{2}}{\partial x\partial z})=0, \\
\frac{\partial f_{3}}{\partial z}+\frac{\partial f_{1}}{\partial y}-x\frac{%
\partial f_{1}}{\partial z}-\lambda \frac{\partial f_{2}}{\partial x}+2(x%
\frac{\partial ^{2}f_{3}}{\partial z\partial x}-\frac{\partial ^{2}f_{3}}{%
\partial x\partial y})=0, \\
\lambda (\lambda \frac{\partial f_{1}}{\partial x}+\frac{\partial f_{2}}{%
\partial y}-x\frac{\partial f_{2}}{\partial z})-2(\frac{\partial ^{2}f_{1}}{%
\partial y\partial z}-x\frac{\partial ^{2}f_{1}}{\partial z^2})=0, \\
\lambda (-\lambda \frac{\partial f_{2}}{\partial x}+\frac{\partial f_{1}}{%
\partial y}-x\frac{\partial f_{1}}{\partial z}+\frac{\partial f_{3}}{%
\partial z})+2(\frac{\partial ^{2}f_{2}}{\partial z\partial y}-x\frac{%
\partial ^{2}f_{2}}{\partial z^2})=0, \\
\lambda (\lambda (f_{1}-\frac{\partial f_{3}}{\partial x})+\frac{\partial
f_{2}}{\partial z})+2(\frac{\partial ^{2}f_{3}}{\partial z\partial y}-x\frac{%
\partial ^{2}f_{3}}{\partial z^2})=0, \\
\frac{\partial f_{1}}{\partial y}-\frac{\partial f_{3}}{\partial z}-x\frac{%
\partial f_{1}}{\partial z}-\lambda \frac{\partial f_{2}}{\partial x}+2\frac{%
\partial ^{2}f_{1}}{\partial z\partial x}=0, \\
\frac{\partial f_{2}}{\partial y}-x\frac{\partial f_{2}}{\partial z}+\lambda
\frac{\partial f_{1}}{\partial x}+2\frac{\partial ^{2}f_{2}}{\partial
z\partial x}=0, \\
\lambda f_{2}+\frac{\partial f_{3}}{\partial y}-x\frac{\partial f_{3}}{%
\partial z}-\frac{\partial f_{1}}{\partial z}+2\frac{\partial ^{2}f_{3}}{%
\partial z\partial x}=0.%
\end{array}
\label{affg2}
\end{equation}%

{By direct calculations, according to the above system of PDEs, we may find $\frac{\partial^2 f_1}{\partial x^2}=\frac{\partial f_1}{\partial y}=\frac{\partial f_1}{\partial z}=0$, so by integrating have}
\begin{equation*}
f_{1}=c_{1}x+c_{2},\text{ }c_{1},c_{2}\in
\mathbb{R}
.
\end{equation*}%
Now, using the thirteenth equation of \eqref{affg2} we obtain, since $\frac{%
\partial f_{1}}{\partial z}=\frac{\partial f_{2}}{\partial z}=0,$ $\frac{%
\partial f_{2}}{\partial y}=-\lambda \frac{\partial f_{1}}{\partial x},$
which gives (after replacing the expression of $f_{1}$)%
\begin{equation*}
f_{2}=-\lambda c_{1}y+c_{3},\text{ }c_{3}\in
\mathbb{R}
.
\end{equation*}%
By replacing $f_{1}$ and $f_{2}$ into the tenth equation of \eqref{affg2} we
prove \ that
\begin{equation*}
f_{3}=\frac{\lambda ^{2}c_{1}}{2}y^{2}-\lambda c_{3}y+A(x),
\end{equation*}%
for a smooth function $A=A(x)$. Hence, the equation {$f_1=\frac{\partial f_3}{\partial x}$} 
implies that $%
A^{\prime }(x)=c_{1}x+c_{2},$ and so,
\begin{equation*}
f_{3}=\frac{\lambda ^{2}c_{1}}{2}y^{2}+\frac{c_{1}}{2}x^{2}+c_{2}x-\lambda
c_{3}y+c_{4},\text{ }c_{4}\in
\mathbb{R}
.
\end{equation*}

So, we proved the following

\begin{theorem}
\label{thg2} $X=f_{1}e_{1}+f_{2}e_{2}+f_{3}e_{3}$ is an affine vector field
on $(H_{3},g_{2})$ if and only if.
\begin{equation}
\left\{
\begin{array}{l}
f_{1}(x,y,z)=c_{1}x+c_{2}, \\
f_{2}(x,y,z)=-\lambda c_{1}y+c_{3},\vphantom{\displaystyle{\frac{A}{B}}} \\
f_{3}(x,y,z)=\frac{\lambda ^{2}c_{1}}{2}y^{2}+\frac{c_{1}}{2}%
x^{2}+c_{2}x-\lambda c_{3}y+c_{4},\vphantom{\displaystyle{\frac{a}{B}}}%
\end{array}%
\right.  \label{affg2comp}
\end{equation}
for any real constants $c_{1},..,c_{4}.$
\end{theorem}

\noindent {Next, one can check easily that the Lie algebra of affine vector
fields of }$(H_{3},g_{2}),$ admits as basis the vector fields given in %
\eqref{g2kill}. Therefore, the following result holds:

\begin{theorem}
\label{thg2a2} Killing vector fields are all and the ones affine vector
fields on $(H_{3},g_{2})$.
\end{theorem}

\subsubsection{{\protect\large \textbf{The third metric}}}

\setcounter{equation}{0}

Next, {let }$X=f_{1}e_{1}+f_{2}e_{2}+f_{3}e_{3}$ be a tangent vector field
on $(H_{3},g_{3})$ with $f_{i}=f_{i}(x,y,z),i=1,2,3$ are arbitrary smooth
functions. We apply \eqref{aff} when $Y,Z$ are chosen in the {orthonormal
basis \eqref{base3}}. Via standard calculations we find that $X$ is an
affine vector field of $(H_{3},g_{3})$ if and only if $f_{1},f_{2},f_{3}$
satisfy the following system of partial differential equations:%
\begin{equation}
\begin{array}{l}
\frac{\partial f_{1}}{\partial x}-2(\frac{\partial f_{2}}{\partial y}+\frac{%
\partial f_{3}}{\partial y})+\frac{\partial ^{2}f_{1}}{\partial y^2}=0, \\
\frac{\partial f_{2}}{\partial x}+2\frac{\partial f_{1}}{\partial y}+\frac{%
\partial ^{2}f_{2}}{\partial y^2}-f_{2}-f_{3}=0, \\
\frac{\partial f_{3}}{\partial x}-2\frac{\partial f_{1}}{\partial y}+\frac{%
\partial ^{2}f_{3}}{\partial y^2}+f_{2}+f_{3}=0, \\
\frac{\partial f_{2}}{\partial x}+\frac{\partial f_{3}}{\partial x}-\frac{%
\partial ^{2}f_{1}}{\partial x\partial y}+\frac{\partial f_{1}}{\partial z}%
=0, \\
\frac{\partial f_{1}}{\partial x}-\frac{\partial f_{2}}{\partial z}+\frac{%
\partial ^{2}f_{2}}{\partial x\partial y}=0, \\
\frac{\partial f_{1}}{\partial x}+\frac{\partial f_{3}}{\partial z}-\frac{%
\partial ^{2}f_{3}}{\partial x\partial y}=0, \\
\frac{\partial ^{2}f_{1}}{\partial y\partial z}-\frac{\partial f_{2}}{%
\partial z}-\frac{\partial f_{3}}{\partial z}=0, \\
\frac{\partial ^{2}f_{2}}{\partial y\partial z}+\frac{\partial f_{1}}{%
\partial z}=0, \\
\frac{\partial ^{2}f_{3}}{\partial y\partial z}-\frac{\partial f_{1}}{%
\partial z}=0, \\
\frac{\partial ^{2}f_{i}}{\partial x^2}=\frac{\partial ^{2}f_{i}}{\partial
z^2}=\frac{\partial ^{2}f_{i}}{\partial x\partial z}=0,\text{ }i=1,2,3.%
\end{array}
\label{affg3}
\end{equation}%
We derive the first equation of \eqref{affg3} with respect to $z,$ and using
the eighth, the ninth and the last equation of \eqref{affg3}, we get $\frac{%
\partial ^{3}f_{1}}{\partial ^{2}y\partial z}=0.$ Integrating, we then have
(since $\frac{\partial ^{2}f_{1}}{\partial ^{2}z}=\frac{\partial ^{2}f_{1}}{%
\partial x\partial z}=0$)

\begin{equation}
\frac{\partial f_{1}}{\partial z}=c_{1}y+c_{2},\text{ }c_{1},c_{2}\in
\mathbb{R}
.  \label{1}
\end{equation}%
Next, taking the derivative of the first and the fourth equation of %
\eqref{affg3} with respect to $x$ and $z,$ respectively, and using the last
equation of \eqref{affg3} and \eqref{1}, we obtain $\frac{\partial ^{3}f_{1}%
}{\partial ^{2}y\partial x}=2c_{1},$ which by integration gives (since $%
\frac{\partial ^{3}f_{1}}{\partial ^{2}y\partial z}=0$)%
\begin{equation}
\frac{\partial ^{2}f_{1}}{\partial ^{2}y}=2c_{1}x+A(y),  \label{2}
\end{equation}%
where $A=A(y)$ is a smooth function.

The first equation of \eqref{affg3} derived with respect to $x,$ gives
(since $\frac{\partial ^{2}f_{1}}{\partial ^{2}x}=0$)%
\begin{equation*}
\frac{\partial ^{2}f_{2}}{\partial x\partial y}+\frac{\partial ^{2}f_{3}}{%
\partial x\partial y}=0.
\end{equation*}%
We then derive the first equation of \eqref{affg3} twice with respect to $y,$
we have, since $\frac{\partial ^{3}f_{1}}{\partial ^{2}y\partial x}=2c_{1},$
$\frac{\partial ^{4}f_{1}}{\partial ^{4}y}=-4c_{1}.$ Thus, from \eqref{2} it
follows%
\begin{equation*}
A(y)=-2c_{1}y^{2}+c_{3}y+c_{4},\text{ }c_{3},c_{4}\in
\mathbb{R}
.
\end{equation*}%
Replacing $A(y)$ into \eqref{2} and integrating we then have
\begin{equation}
\frac{\partial f_{1}}{\partial y}=(2c_{1}x+c_{4})y-\frac{2c_{1}}{3}y^{3}+%
\frac{c_{3}}{2}y^{2}+B(x,z),  \label{3}
\end{equation}%
for a smooth function $B=B(x,z).$

Next derive \eqref{1} with respect to $z$ and taking into account \eqref{3},
we prove that $B=c_{1}z+C(x),$ for a smooth function $C=C(x).$ We then
derive \eqref{3} with respect to $x,$ we get (since $\frac{\partial ^{2}f_{1}%
}{\partial ^{2}x}=0$) $C^{\prime \prime }(x)=0,$ that is $C(x)=c_{5}x+c_{6},$
$c_{5},c_{6}\in
\mathbb{R}
.$ Therefore equation \eqref{3} becomes%
\begin{equation*}
\frac{\partial f_{1}}{\partial y}=(2c_{1}x+c_{4})y-\frac{2c_{1}}{3}y^{3}+%
\frac{c_{3}}{2}y^{2}+c_{1}z+c_{5}x+c_{6}.
\end{equation*}%
Integrating, we find%
\begin{equation}
f_{1}=-\frac{c_{1}}{6}y^{4}+\frac{c_{3}}{6}y^{3}+\frac{(2c_{1}x+c_{4})}{2}%
y^{2}+\left( c_{1}z+c_{5}x+c_{6}\right) y+D(x,z),  \label{4}
\end{equation}%
for a smooth function $D=D(x,z).$ From \eqref{1}, \eqref{4} and using the
fact that $\frac{\partial ^{2}f_{1}}{\partial x\partial z}=\frac{\partial
^{2}f_{1}}{\partial ^{2}z}=\frac{\partial ^{2}f_{1}}{\partial ^{2}x}=0,$ we
easily prove that $D=c_{2}z+c_{7}x+c_{8},c_{7},c_{8}$ are real constants.
Hence, \eqref{4} becomes%
\begin{equation}
f_{1}=-\frac{c_{1}}{6}y^{4}+\frac{c_{3}}{6}y^{3}+\frac{(2c_{1}x+c_{4})}{2}%
y^{2}+\left( c_{1}z+c_{5}x+c_{6}\right) y+c_{2}z+c_{7}x+c_{8}.  \label{5}
\end{equation}%
Now, replacing \eqref{5} into the eighth equation of \eqref{affg3} and
taking into account that $\frac{\partial ^{2}f_{2}}{\partial x\partial z}=%
\frac{\partial ^{2}f_{2}}{\partial ^{2}z}=\frac{\partial ^{2}f_{2}}{\partial
^{2}x}=0$, we easily get%
\begin{equation*}
\frac{\partial f_{2}}{\partial z}=-\frac{c_{1}}{2}y^{2}-c_{2}y+c_{9},c_{9}%
\in
\mathbb{R}
,
\end{equation*}%
which by integration yields%
\begin{equation}
f_{2}=\left( -\frac{c_{1}}{2}y^{2}-c_{2}y+c_{9}\right) z+G(y)x+H(y),
\label{6}
\end{equation}%
for arbitrary real smooth functions $G=G(y)$ and $H=H(y).$

By replacing \eqref{5} and \eqref{6} into the fifth equation of \eqref{affg3}%
, we obtain%
\begin{equation*}
G(y)=-\frac{c_{1}}{2}y^{3}-\frac{(c_{5}+c_{2})}{2}y^{2}+\left(
c_{9}-c_{7}\right) y+c_{10},\text{ }c_{10}\in
\mathbb{R}
.
\end{equation*}%
Taking the derivative of the second equation of \eqref{affg3} with respect
to $y$. By the first equation of \eqref{affg3} and using the expression of $%
f_{1},f_{2}$ and $G,$ we prove that
\begin{equation*}
H^{\prime \prime \prime }(y)=5c_{1}y^{2}+\left( \frac{3c_{5}-3c_{3}+2c_{2}}{2%
}\right) y+\frac{3c_{7}-3c_{4}-2c_{9}}{2}.
\end{equation*}%
Integrating, we then prove that equation \eqref{6} becomes
\begin{eqnarray}
f_{2} &=&\left( -\frac{c_{1}}{2}y^{2}-c_{2}y+c_{9}\right) z+\left( -\frac{%
c_{1}}{2}y^{3}-\frac{(c_{5}+c_{2})}{2}y^{2}+(c_{9}-c_{7})y+c_{10}\right) x+%
\frac{c_{1}}{12}y^{5}  \label{7} \\
&&+\left( \frac{3c_{5}-3c_{3}+2c_{2}}{48}\right) y^{4}+\left( \frac{%
3c_{7}-3c_{4}-2c_{9}}{12}\right) y^{3}+\frac{c_{11}}{2}y^{2}+c_{12}y+c_{13},
\notag
\end{eqnarray}%
with $c_{11},c_{12},c_{13}\in
\mathbb{R}
.$ Therefore, the second equation of \eqref{affg3} together with \eqref{5}
and \eqref{7}, yields
\begin{eqnarray}
f_{3} &=&\left( \frac{c_{1}}{2}y^{2}+c_{2}y+c_{1}-c_{9}\right) z-\frac{c_{1}%
}{12}y^{5}-\left( \frac{3c_{5}-3c_{3}+2c_{2}}{48}\right) y^{4}  \label{8} \\
&&+\left( \frac{-3c_{7}+3c_{4}+2c_{9}-2c_{1}}{12}\right) y^{3}+\left( \frac{%
c_{3}+c_{5}-2c_{11}}{4}\right) y^{2}  \notag \\
&&+\left( \frac{c_{1}}{2}y^{3}+\frac{(c_{2}+c_{5})}{2}%
y^{2}+(c_{1}+c_{7}-c_{9})y+c_{5}-c_{2}-c_{10}\right) x  \notag \\
&&+\left( \frac{c_{4}+c_{7}-2c_{12}}{2}\right) y+c_{10}+c_{11}+2c_{6}-c_{13}.
\notag
\end{eqnarray}%
Replacing \eqref{5}, \eqref{7} and \eqref{8} into the third equation of %
\eqref{affg3}. Standard calculations show that $c_{3}=2c_{2}-3c_{5}.$ This
leads to prove the following.

\begin{theorem}
$X=f_{1}e_{1}+f_{2}e_{2}+f_{3}e_{3}$ is an affine vector field on $%
(H_{3},g_{3})$ if and only if.
\begin{equation}
\left\{
\begin{array}{l}
f_{1}(x,y,z)=-\frac{c_{1}}{6}y^{4}+\frac{\left( 2c_{2}-3c_{5}\right) }{6}%
y^{3}+\frac{(2c_{1}x+c_{4})}{2}y^{2}+\left( c_{1}z+c_{5}x+c_{6}\right)
y+c_{2}z+c_{7}x+c_{8}, \\
\\
f_{2}(x,y,z)=\left( -\frac{c_{1}}{2}y^{2}-c_{2}y+c_{9}\right) z+\left( -%
\frac{c_{1}}{2}y^{3}-\frac{(c_{5}+c_{2})}{2}y^{2}+(c_{9}-c_{7})y+c_{10}%
\right) x+\frac{c_{1}}{12}y^{5} \\
\text{ \ \ \ \ \ \ \ \ \ \ \ \ \ \ \ \ \ \ }+\left( \frac{3c_{5}-c_{2}}{12}%
\right) y^{4}+\left( \frac{3c_{7}-3c_{4}-2c_{9}}{12}\right) y^{3}+\frac{%
c_{11}}{2}y^{2}+c_{12}y+c_{13}, \\
\\
f_{3}(x,y,z)=\left( \frac{c_{1}}{2}y^{2}+c_{2}y+c_{1}-c_{9}\right) z-\frac{%
c_{1}}{12}y^{5}-\left( \frac{3c_{5}-c_{2}}{12}\right) y^{4}+\left( \frac{%
-3c_{7}+3c_{4}+2c_{9}-2c_{1}}{12}\right) y^{3} \\
\text{ \ \ \ \ \ \ \ \ \ \ \ \ \ \ \ \ \ \ }+\left( \frac{c_{1}}{2}y^{3}+%
\frac{(c_{2}+c_{5})}{2}y^{2}+(c_{1}+c_{7}-c_{9})y+c_{5}-c_{2}-c_{10}\right)
x+\left( \frac{c_{2}-c_{5}-c_{11}}{2}\right) y^{2} \\
\text{ \ \ \ \ \ \ \ \ \ \ \ \ \ \ \ \ \ \ }+\left( \frac{c_{4}+c_{7}-2c_{12}%
}{2}\right) y+c_{10}+c_{11}+2c_{6}-c_{13},%
\end{array}%
\right.  \label{affg3comp}
\end{equation}%
for any real constants $c_{i},i=1,..13$ and $i\neq 3.$
\end{theorem}

\section{{\protect\large \textbf{Ricci, curvature and matter collineations }}%
}

In this Section, we shall investigate symmetries of the three-dimensional
Heisenberg group $H_{3}$, equipped with each of left-invariant metrics $g_0$, $g_1$, $g_2$.
%
collineation.

\subsection{\protect\large Riemannian setting}

In this case, by using to the orthonormal basis $\{e_1,e_2,e_3\}$, we
firstly calculate the Lie derivative of the Ricci tensor with respect to the
smooth vector field $X=f_1e_1+\cdots+f_3e_3$, that is
\begin{equation*}
\begin{tabular}{l}
$\left(\mathcal{L}_X\varrho\right)\left(e_1,e_1\right)=\lambda^2\left(x\frac{%
\partial f_1}{\partial z}-\frac{\partial f_1}{\partial y}\right)$, \\
$\left(\mathcal{L}_X\varrho\right)\left(e_1,e_2\right)=\frac{\lambda^2}{2}%
\left(x\frac{\partial f_2}{\partial z}-\lambda\frac{\partial f_1}{\partial x}%
-\frac{\partial f_2}{\partial y}\right)$, \\
$\left(\mathcal{L}_X\varrho\right)\left(e_1,e_3\right)=-\frac{\lambda^2}{2}%
\left(x\frac{\partial f_3}{\partial z}-\frac{\partial f_3}{\partial y}+\frac{%
\partial f_1}{\partial z}-\lambda f_2\right)$, \\
$\left(\mathcal{L}_X\varrho\right)\left(e_2,e_2\right)=-\lambda^3\frac{%
\partial f_2}{\partial x}$, \\
$\left(\mathcal{L}_X\varrho\right)\left(e_2,e_3\right)=\frac{\lambda^2}{2}%
\left(\lambda\frac{\partial f_3}{\partial x}-\lambda f_1-\frac{\partial f_2}{%
\partial z}\right)$, \\
$\left(\mathcal{L}_X\varrho\right)\left(e_3,e_3\right)=\lambda^2\frac{%
\partial f_3}{\partial z}$. \\
\end{tabular}%
\end{equation*}
By the above relations, $X$ is a curvature collineation if and only if be a
Ricci collineation if and only if
\begin{equation}  \label{Ricog0}
\left\{%
\begin{array}{l}
f_1(x,y,z)=-\frac{c_1}{\lambda}x+c_3, \\
f_2(x,y,z)=c_1y+c_2, \\
f_3(x,y,z)=-\frac{c_1\lambda}{2}y^2-\lambda c_2y-\frac{c_1}{2\lambda}%
x^2+c_3x+c_4,%
\end{array}%
\right.
\end{equation}
for arbitrary real constants $c_1,\cdots c_4$. To study matter
collineations, by using the equations \eqref{g0}, \eqref{ro0} and \eqref{to0}%
, we firstly determine the tensor field $T$ as following
\begin{equation*}
T=\left(%
\begin{array}{ccc}
-\frac{\lambda^2}{4} & 0 & 0 \\
0 & -\frac{\lambda^2}{4} & 0 \\
0 & 0 & \frac{3\lambda^2}{4}%
\end{array}
\right).
\end{equation*}
We then calculate the Lie derivatives of the above tensor field $T$ as
\begin{equation*}
\begin{array}{l}
\left(\mathcal{L}_XT\right)(e_1,e_1)=\frac{\lambda^2}{2}\left(x\frac{%
\partial f_1}{\partial x}-\frac{\partial f_1}{\partial y}\right), \\
\left(\mathcal{L}_XT\right)=-\frac{\lambda^2}{4}\left(\frac{\partial f_2}{%
\partial y}-x\frac{\partial f_2}{\partial z}+\lambda\frac{\partial f_1}{%
\partial x}\right), \\
\left(\mathcal{L}_XT\right)=-\frac{\lambda^2}{4}\left(3(x\frac{\partial f_3}{%
\partial z}-\frac{\partial f_3}{\partial y}-\lambda f_2)+\frac{\partial f_1}{%
\partial z}\right), \\
\left(\mathcal{L}_XT\right)=-\frac{\lambda^3}{2}\frac{\partial f_2}{\partial
x}, \\
\left(\mathcal{L}_XT\right)=-\frac{\lambda^2}{4}\left(3\lambda(f_1-\frac{%
\partial f_3}{\partial x})+\frac{\partial f_2}{\partial z}\right), \\
\left(\mathcal{L}_XT\right)=\frac{3\lambda^2}{2}\frac{\partial f_3}{\partial
z}.%
\end{array}%
\end{equation*}
By the above equations Matter collineations will coincide with the Ricci
collineations, i.e., comparing with the results included in the previous Section, $X$ is a Killing vector fields and so we don't have any
nontrivial matter collineation in this case.

\subsection{\protect\large Lorentzian setting}

\subsubsection{\protect\large The first metric}

With respect to the orthonormal basis $\{e_{1},e_{2},e_{3}\}$ given by %
\eqref{base1}, l{et }$X=f_{1}e_{1}+f_{2}e_{2}+f_{3}e_{3}$ denote an
arbitrary vector field on $(H_{3},g_{1}),$ where $f_{i}=f_{i}(x,y,z),i=1,2,3$
are arbitrary smooth functions. We now determine the Lie derivative $\left(
\mathcal{L}_{X}\varrho \right) \left( e_{i},e_{j}\right) ,$ for all indices $%
i\leq j.$ We have the following%
\begin{equation*}
\begin{tabular}{l}
$\left( \mathcal{L}_{X}\varrho \right) \left( e_{1},e_{1}\right) =-\lambda
^{2}\left( \frac{\partial f_{1}}{\partial z}\right) ,$ \\
$\left( \mathcal{L}_{X}\varrho \right) \left( e_{1},e_{2}\right) =\frac{%
\lambda ^{2}}{2}\left( \frac{\partial f_{2}}{\partial z}-\lambda f_{3}-\frac{%
\partial f_{1}}{\partial y}+x\frac{\partial f_{1}}{\partial z}\right) ,$ \\
$\left( \mathcal{L}_{X}\varrho \right) \left( e_{1},e_{3}\right) =\frac{%
\lambda ^{2}}{2}\left( -\frac{\partial f_{3}}{\partial z}+\lambda
f_{2}-\lambda \frac{\partial f_{1}}{\partial x}\right) ,$ \\
$\left( \mathcal{L}_{X}\varrho \right) \left( e_{2},e_{2}\right) =\lambda
^{2}\left( \frac{\partial f_{2}}{\partial y}-x\frac{\partial f_{2}}{\partial
z}\right) ,$ \\
$\left( \mathcal{L}_{X}\varrho \right) \left( e_{2},e_{3}\right) =\frac{%
\lambda ^{2}}{2}\left( -\frac{\partial f_{3}}{\partial y}+x\frac{\partial
f_{3}}{\partial z}+\lambda \frac{\partial f_{2}}{\partial x}\right) ,$ \\
$\left( \mathcal{L}_{X}\varrho \right) \left( e_{3},e_{3}\right) =-\lambda
^{3}\left( \frac{\partial f_{3}}{\partial x}\right) .$%
\end{tabular}%
\end{equation*}%
Ricci collineations are then calculated by solving the systeme of PDE
obtained by requiring that all the above components of $\mathcal{L}%
_{X}\varrho $ vanish. This proves the following result:

\begin{theorem}
Let $X=f_{1}e_{1}+f_{2}e_{2}+f_{3}e_{3}$ denote an arbitrary vector field on
the Heisenberg group $(H_{3},g_{1}).$ Then, $X$ is a Ricci collineation if
and only if
\begin{equation}
\left\{
\begin{array}{l}
f_{1}(x,y,z)=\frac{c_{1}}{2\lambda }x^{2}-\frac{\lambda c_{1}}{2}%
y^{2}+c_{2}x-\lambda c_{3}y+c_{4}, \\
f_{2}(x,y,z)=\frac{c_{1}}{\lambda }x+c_{2}, \\
f_{3}(x,y,z)=c_{1}y+c_{3},%
\end{array}%
\right.  \label{Ricog1}
\end{equation}%
for any real constants $c_{1},c_{2},c_{3},c_{4}$.
\end{theorem}

\begin{remark}
If we set $X=f_{1}e_{1}+f_{2}e_{2}+f_{3}e_{3},$ where $f_{1},f_{2},f_{3}$
are smooth functions given by \eqref{Ricog1}, From {\eqref{g1kill} }we have
that $X$ is a Killing vector field of $(H_{3},g_{1}).$ Therefore, since we
are in dimension $3$, Ricci collineation of $(H_{3},g_{1})$ will coincide
with the curvature collineation.
\end{remark}

We shall now classify matter collineations on the Heisenberg
group $(H_{3},g_{1}).$ Starting from \eqref{g1}, \eqref{ro1} and \eqref{to1}%
, a direct calculations yields that for the Lorentzian metric $g_{1},$ with
respect to the orthonormal basis given in \eqref{base1} the tensor field $T$
is described by%
\begin{equation*}
T=\left(
\begin{array}{ccc}
-\frac{3\lambda ^{2}}{4} & 0 & 0 \\
0 & \frac{\lambda ^{2}}{4} & 0 \\
0 & 0 & -\frac{\lambda ^{2}}{4}%
\end{array}%
\right) .
\end{equation*}%
Next, let $X=f_{1}e_{1}+f_{2}e_{2}+f_{3}e_{3}$ denote an arbitrary vector
field on $(H_{3},g_{1}),$ where $f_{i}=f_{i}(x,y,z),i=1,2,3$ are arbitrary
smooth functions. We then compute the Lie derivative of $T$ with respect to $%
X$ and we find%
\begin{equation*}
\begin{tabular}{l}
$\left( \mathcal{L}_{X}T\right) \left( e_{1},e_{1}\right) =-\frac{3\lambda
^{2}}{2}\left( \frac{\partial f_{1}}{\partial z}\right) ,$ \\
$\left( \mathcal{L}_{X}T\right) \left( e_{1},e_{2}\right) =\frac{\lambda ^{2}%
}{4}\left( \frac{\partial f_{2}}{\partial z}-3\left( \lambda f_{3}+\frac{%
\partial f_{1}}{\partial y}-x\frac{\partial f_{1}}{\partial z}\right)
\right) ,$ \\
$\left( \mathcal{L}_{X}T\right) \left( e_{1},e_{3}\right) =\frac{\lambda ^{2}%
}{4}\left( -\frac{\partial f_{3}}{\partial z}+3\lambda \left( f_{2}-\frac{%
\partial f_{1}}{\partial x}\right) \right) ,$ \\
$\left( \mathcal{L}_{X}T\right) \left( e_{2},e_{2}\right) =\frac{\lambda ^{2}%
}{2}\left( \frac{\partial f_{2}}{\partial y}-x\frac{\partial f_{2}}{\partial
z}\right) ,$ \\
$\left( \mathcal{L}_{X}T\right) \left( e_{2},e_{3}\right) =\frac{\lambda ^{2}%
}{4}\left( -\frac{\partial f_{3}}{\partial y}+x\frac{\partial f_{3}}{%
\partial z}+\lambda \frac{\partial f_{2}}{\partial x}\right) ,$ \\
$\left( \mathcal{L}_{X}T\right) \left( e_{3},e_{3}\right) =-\frac{\lambda
^{3}}{2}\left( \frac{\partial f_{3}}{\partial x}\right) .$%
\end{tabular}%
\end{equation*}%
Requiring that all above components of $\mathcal{L}_{X}T$ vanish, we have
that, $X$ is a matter collineation if and only if $X$ is a Ricci
collineation (that is, $X$ a Killing vector field). Therefore, there is no
nontrivial matter collineation for $(H_{3},g_{1}).$

\subsubsection{\protect\large The second metric}

We now determine Ricci collineation of $(H_{3},g_{2})$. Put $%
X=f_{1}e_{1}+f_{2}e_{2}+f_{3}e_{3},$ where $f_{i}=f_{i}(x,y,z),i=1,2,3$ are
arbitrary smooth functions, and $\{e_{1},e_{2},e_{3}\}$ is the orthonormal
basis given by \eqref{base2}. We now determine the Lie derivative $\left(
\mathcal{L}_{X}\varrho \right) \left( e_{i},e_{j}\right) ,$ for all indices $%
i\leq j.$ We have the following%
\begin{equation*}
\begin{tabular}{l}
$\left( \mathcal{L}_{X}\varrho \right) \left( e_{1},e_{1}\right) =\lambda
^{2}\left( \frac{\partial f_{1}}{\partial y}-x\frac{\partial f_{1}}{\partial
z}\right) ,$ \\
$\left( \mathcal{L}_{X}\varrho \right) \left( e_{1},e_{2}\right) =\frac{%
\lambda ^{2}}{2}\left( \frac{\partial f_{2}}{\partial y}-x\frac{\partial
f_{2}}{\partial z}+\lambda \frac{\partial f_{1}}{\partial x}\right) ,$ \\
$\left( \mathcal{L}_{X}\varrho \right) \left( e_{1},e_{3}\right) =\frac{%
\lambda ^{2}}{2}\left( \lambda f_{2}+\frac{\partial f_{3}}{\partial y}+\frac{%
\partial f_{1}}{\partial z}-x\frac{\partial f_{3}}{\partial z}\right) ,$ \\
$\left( \mathcal{L}_{X}\varrho \right) \left( e_{2},e_{2}\right) =\lambda
^{3}\frac{\partial f_{2}}{\partial x},$ \\
$\left( \mathcal{L}_{X}\varrho \right) \left( e_{2},e_{3}\right) =\frac{%
\lambda ^{2}}{2}\left( \frac{\partial f_{2}}{\partial z}+\lambda \frac{%
\partial f_{3}}{\partial x}-\lambda f_{1}\right) ,$ \\
$\left( \mathcal{L}_{X}\varrho \right) \left( e_{3},e_{3}\right) =\lambda
^{2}\frac{\partial f_{3}}{\partial z}.$%
\end{tabular}%
\end{equation*}%
Ricci collineations are then calculated by solving the systeme of PDE
obtained by requiring that all the above coefficients of $\mathcal{L}%
_{X}\varrho $ vanish. This leads to prove the following :

\begin{theorem}
Let $X=f_{1}e_{1}+f_{2}e_{2}+f_{3}e_{3}$ denote an arbitrary vector field on
the Heisenberg group $(H_{3},g_{2}).$ Then, $X$ is a Ricci collineation if
and only if
\begin{equation}
\left\{
\begin{array}{l}
f_{1}(x,y,z)=-\frac{c_{1}}{\lambda }x+c_{2}, \\
f_{2}(x,y,z)=c_{1}y+c_{3}, \\
f_{3}(x,y,z)=-\frac{c_{1}}{2\lambda }x^{2}-\frac{\lambda c_{1}}{2}%
y^{2}+c_{2}x-\lambda c_{3}y+c_{4},%
\end{array}%
\right.  \label{Ricog2}
\end{equation}%
for arbitrary real constants $c_{1},c_{2},c_{3},c_{4}$.
\end{theorem}

\begin{remark}
Put $X=f_{1}e_{1}+f_{2}e_{2}+f_{3}e_{3},$ where $f_{1},f_{2},f_{3}$ are
smooth functions given by \eqref{Ricog2}, From {\eqref{g2kill} it follows }%
that $X$ is a Killing vector field of $(H_{3},g_{2}).$ Thus, from the fact
that $\dim H_{3}=3,$ we prove that $X$ is a Ricci collineation of $%
(H_{3},g_{2})$ if and only if $X$ is a curvature collineation.
\end{remark}

Now, consider the Heisenberg group endowed with the Lorentzian metric $%
g_{2}. $ From \eqref{g2}, \eqref{ro2} and \eqref{to2}, it follows that the
tensor field $T$ is given, with respect to the orthonormal basis %
\eqref{base2}, by%
\begin{equation*}
T=\left(
\begin{array}{ccc}
\frac{\lambda ^{2}}{4} & 0 & 0 \\
0 & \frac{\lambda ^{2}}{4} & 0 \\
0 & 0 & \frac{3\lambda ^{2}}{4}%
\end{array}%
\right) .
\end{equation*}%
Next, put $X=f_{1}e_{1}+f_{2}e_{2}+f_{3}e_{3},$ where $%
f_{i}=f_{i}(x,y,z),i=1,2,3$ are arbitrary smooth functions. We then compute
the Lie derivative of $T$ with respect to $X$ and we get%
\begin{equation*}
\begin{tabular}{l}
$\left( \mathcal{L}_{X}T\right) \left( e_{1},e_{1}\right) =\frac{\lambda ^{2}%
}{2}\left( \frac{\partial f_{1}}{\partial y}-x\frac{\partial f_{1}}{\partial
z}\right) ,$ \\
$\left( \mathcal{L}_{X}T\right) \left( e_{1},e_{2}\right) =\frac{\lambda ^{2}%
}{4}\left( \frac{\partial f_{2}}{\partial y}-x\frac{\partial f_{2}}{\partial
z}+\lambda \frac{\partial f_{1}}{\partial x}\right) ,$ \\
$\left( \mathcal{L}_{X}T\right) \left( e_{1},e_{3}\right) =\frac{\lambda ^{2}%
}{4}\left( \frac{\partial f_{1}}{\partial z}+3\left( \lambda f_{2}+\frac{%
\partial f_{3}}{\partial y}-x\frac{\partial f_{3}}{\partial z}\right)
\right) ,$ \\
$\left( \mathcal{L}_{X}T\right) \left( e_{2},e_{2}\right) =\frac{\lambda ^{3}%
}{2}\left( \frac{\partial f_{2}}{\partial x}\right) ,$ \\
$\left( \mathcal{L}_{X}T\right) \left( e_{2},e_{3}\right) =\frac{\lambda ^{2}%
}{4}\left( \frac{\partial f_{2}}{\partial z}-3\lambda \left( f_{1}+\frac{%
\partial f_{3}}{\partial x}\right) \right) ,$ \\
$\left( \mathcal{L}_{X}T\right) \left( e_{3},e_{3}\right) =\frac{3\lambda
^{2}}{2}\left( \frac{\partial f_{3}}{\partial z}\right) .$%
\end{tabular}%
\end{equation*}%
Solving the system of PDEs obtained by requiring that $\mathcal{L}_{X}T=0,$
we get that all solutions coincide with the Killing vector fields of $%
(H_{3},g_{2}).$ This means that matter collineations are trivial.

\subsubsection{\protect\large The third metric}

In this case, since the metric $g_{3}$ is flat, we deduce that all vector
fields of $(H_{3},g_{3})$ are Ricci, and curvature collineation. Moreover,
matter collineations coincide with Killing vector fields.

\end{document}